\documentclass[12pt]{article}

\usepackage{amsfonts}
\usepackage{amssymb}
\usepackage{amsmath}
\usepackage{amsthm}

\def \({\left(}
\def \){\right)}
\def \le {\leqslant}
\def \ge {\geqslant}

\begin{document}
\begin{large}
\centerline{\bf On certain Littlewood-like and Schmidt-like problems}
\centerline{{\bf in inhomogeneous Diophantine approximations} }
\end{large}
\vskip+1.0cm

\centerline{ by {\bf Nikolay  Moshchevitin}\footnote{supported by the grant RFBR № 09-01-00371
  }}
\vskip+1.0cm

We give several results related to inhomogeneous approximations  to two real numbers and
badly approximable numbers. Our results are related to classical theorems by A. Khintchine \cite{k26}
and to an original method invented by Y. Peres and W. Schlag \cite{PS}.

\vskip+1.0cm
 {\bf 1.  Functions and parameters.}

In all what follows, $||\cdot ||$ is the distance to the nearest integer.
 All functions
here are non-negative valued functions in real non-negative
variables.

Consider strictly increasing functions $\omega_1(t), \omega_2
(t)$. Let $\omega_1^*(t)$ be the inverse function to
$\omega_1(t)$, that is
$$
\omega_1^*(\omega_1 (t)) = t. 
$$
Suppose that another function in two variables function
$\Omega(x,y)$ satisfy the condition
\begin{equation}\label{amega}
 \begin{cases}
 x y \le \omega \left(\frac{z}{x}\right),\cr x\le z
\end{cases}
 \,\,\,\,\,\, \Longrightarrow\,\,\,\,\,\,
x\le \Omega( y, z), \,\,\, \,\,\,\forall \, x,y,z \in
\mathbb{Z}_+.
\end{equation}
This condition may be rewritten as
\begin{equation}\label{amega00}
 \begin{cases}
 x\omega_1^*
(x\cdot y) \le z,\cr x\le z
\end{cases}
 \,\,\,\,\,\, \Longrightarrow\,\,\,\,\,\,
x\le \Omega( y, z), \,\,\, \,\,\,\forall \, x,y,z \in
\mathbb{Z}_+.
\end{equation}

 Suppose that the
functions $  \phi (t), \phi_2 (t), \phi_2 (t), \psi_1(t),
\psi_2(t),$ increase as $t \to \infty$ and
\begin{equation}\label{zero}
\phi(0) = \phi_1(0)= \phi_2(0) = \psi_1(0)= \psi_2(0)= 0 .
\end{equation}

Suppose that $\psi_j(t),\,  j= 1,2$ are strictly increasing
functions  and that $\psi_j^* (t)$ is the inverse function of
$\psi_j(t)$, that is
$$
\psi_j^* (\psi_j(t)) = t\,\,\,\, \forall \, t \in
\mathbb{R}_+,\,\,\,\,\, j = 1,2.
$$

For a positive $\varepsilon >0$ and  integers $\nu, \mu$ define  
 \begin{equation}\label{delta1}
   \delta^{[1]}_\varepsilon (\mu,\nu) =
\psi_2^*\left(  \frac{ \varepsilon}{\phi (2^\nu) \psi_1
(2^{-\mu-1})}
 \right),
\end{equation}
 \begin{equation}\label{delta12}
   \delta^{[2]}_\varepsilon (\nu) =
\psi_2^*\left(  \frac{ \varepsilon}{\phi_2 (2^\nu)  }
 \right).
\end{equation}

Suppose that $A>1$.  For functions $\omega_1 (t),\omega_2(t),
\phi(t), \psi_1(t) , \psi_2(t)$  we consider the following sum:

\begin{equation}\label{docond}
 S_{A,\varepsilon}^{[1]}(X)=\,\,\,\sum_{X\le \nu <A(X+1)}\, \sum_{1\le \mu \le
\log_2 (\omega_2 (2^{\nu+1}))+1}\, \delta^{[1]}_\varepsilon
(\mu,\nu)\cdot
  \max\left(\Omega (2^{\mu-1},
2^{\nu +1}), 2^{\nu  -\mu},1 \right)
\end{equation}

For functions $\omega_1 (t),\omega_2(t) , \phi_1(t), \phi_2(t),  \psi_1(t) , \psi_2(t)$ we
consider another sum:

\begin{equation}\label{docond2}
 S_{A,\varepsilon}^{[2]}(X)=\,\,\,\sum_{X\le \nu <A(X+1)}\,\,\,
\delta^{[2]}_\varepsilon
(\nu)\cdot
\max\left(\Omega
(1/2r_\varepsilon (\nu ),2^{\nu+1}) ,2^{\nu}r_\varepsilon (\nu ) , 1\right),
\end{equation}
where

\begin{equation}\label{docond3}
r_\varepsilon (\nu ) =
\psi_1^*\left(\frac{\varepsilon}{\phi_1(2^{\nu})}\right).
\end{equation}

{\bf 2. Main results.}

Here we formulate two new results - Theorems 1,2.
 Proofs of this theorems are given in Sections 6, 7, 8.
 Section 4 below is devoted to certain examples of applications of
 Theorem 1.
 Section 5 deals with applications of Theorem 2.
 In Section 3, we discuss Khintchine's theorems and some
 of their extensions.

 {\bf Theorem 1.}

{\it  Suppose that functions $ \psi_1(t), \psi_2(t), \phi (t)
 $ are increasing. Suppose that (\ref{zero}) is valid.
  Suppose that for certain $ A>1,\varepsilon
>0, X_0 \ge 0$  all the functions satisfy the conditions
\begin{equation}\label{cond0}
  \log_2 \left(\frac{ X}{2\psi_2^*\left(\frac{\varepsilon}{\phi(X)\psi_1(1/2)}\right)}\right)\le
  (A-1)\log_2 X,\,\,\,\,\forall X\ge X_0,
\end{equation}
 and
\begin{equation}\label{cond}
 \sup_{X\ge X_0} \,\,\, S^{[1]}_{A,\varepsilon}(X)\le
\frac{1}{2^9}.
\end{equation}

Consider two real numbers $\alpha, \eta$ such that
\begin{equation}\label{bad}
\inf_{x\ge X_0}\, \omega_1 (x)\cdot ||x\alpha || \ge 1
\end{equation}
and
\begin{equation}\label{bad1}
\inf_{x\ge X_0}\, \omega_2 (x)\cdot ||x\alpha -\eta || \ge 1
\end{equation}

 Then for any sequence of real numbers $
\eta_1,\eta_2,...,\eta_x,... $ there exists a real number $\beta$
such that
\begin{equation}\label{finale}
\inf_{x\ge X_0}\,  \phi (x) \psi_1(||x\alpha -\eta||)
 \psi_2(||x\beta -\eta_x||)\,  \ge \varepsilon.
\end{equation}
}

 A simpler version of the theorem was announced in \cite{MB} (Theorem 8 from \cite{MB}).
 Some inhomogeneous results in special case were announced in
\cite{U} (see Appendix from \cite{U}).

The following  Theorem 2 generalizes a result from \cite{ME}.

{\bf Theorem 2.} {\it

Consider a real number $\alpha$  satisfying
(\ref{bad}). Let $\eta $ be an arbitrary real number.
Suppose that
\begin{equation}\label{cond000}
  \log_2 \left(\frac{ X}{2\psi_2^*\left(\frac{\varepsilon}{\phi_2(X)}\right)}\right)\le
  (A-1)\log_2 X,\,\,\,\,\forall X\ge X_0,
\end{equation}
 and
\begin{equation}\label{condmore}
 \sup_{X\ge X_0} \,\,\, S^{[2]}_{A,\varepsilon}(X)\le
\frac{1}{2^9}.
\end{equation}

 Then for any sequence of real numbers $
\eta_1,\eta_2,...,\eta_x,... $ there exists a real number $\beta$
such that
\begin{equation}\label{finale1}
\inf_{x\ge X_0}\, \max \left( \phi_1 (x)\cdot \psi_1(||x\alpha
-\eta||),\,\,
 \phi_2 (x)\cdot\psi_2(||x\beta -\eta_x||)\,\right)\,  \ge \varepsilon.
\end{equation}
}

 {\bf Remark.} The method under consideration enables one to obtain
results about intersections. Suppose that $ j \in \{1, 2\}$.
 Given two different
collections of functions $\omega_1^j (t), \omega_2^j(t),\psi_1^j
(t), \psi_2^j(t) ,\phi^j (t),  \sigma_1^j(t) ,\sigma_2^j (t)$, two
sequences $\{\eta_x^j\}_{x=1}^\infty$ and two couples of reals
$\alpha^j,\eta^j$ satisfying the conditions specified (with more
restrictions on constants) it is easy to prove  the existence of a
real $\beta$ such that the conclusions (\ref{finale},
\ref{finale1}) (or even both of them) are valid for  both 
values of $ j \in \{1, 2\}$. A simpler example of such a result  was
proved in \cite{ME}. Moreover the method can give lower bound for
Hausdorff dimension of the sets.

 {\bf 3.
Khintchine's theorems and their extensions.}

In \cite{k26} A. Khintchine proved the following result.

{\bf Theorem A.} {\it There exists an absolute constant $\gamma$
such that for any real $\alpha$ there exists a real $\eta$ such
that
\begin{equation}\label{kin}
\inf_{x\in \mathbb{Z}_+}\,  x\cdot ||x\alpha -\eta || \ge \gamma .
\end{equation}
 }

One can find this theorem in the books \cite{Kas} (Ch. 10) and
\cite{SUZ} (Ch. 4). The best known value of $\gamma $ probably is
due to H. Godwin \cite{god}. From \cite{Ts} we know that for every
$\alpha \in \mathbb{R}$ the set of all $\eta$ for which there
exists a positive constant $\gamma$ such that (\ref{kin}) is true
is a 1/2-winning set.

From Khintchine's theorem it follows that there exist reals
$\alpha, \eta$ such that inequalities (\ref{bad}),  (\ref{bad1}) are
valid with
$$
\omega_1 (t) = \omega_2 (t) = \gamma t
$$
with an absolute positive constant $\gamma$.

Here we formulate an immediate corollary to Khintchine's Theorem
A.

{\bf Corollary 1.}\,\, {\it

{\rm (i)} Suppose that  reals $\alpha_1$  and $\alpha_2$ are
linearly dependent over $\mathbb Z$ together with $1$. Them there
exist reals $\eta_1, \eta_2$ such that
$$
\inf_{x\in \mathbb{Z}_+}\,  x\cdot ||x\alpha_1 -\eta_1 || \cdot
||x\alpha_2 -\eta_2 ||
>0.
$$

{\rm (ii)} Suppose that $\alpha_1$ is a badly approximable number
satisfying
$$
\inf_{x\in \mathbb{Z}_+}\,  x\cdot ||x\alpha_1 ||
>0.
$$
Suppose that $\alpha_2$ is linearly dependent with $\alpha_1$ and
$1$. Then there exists  $ \eta$ such that
$$
\inf_{x\in \mathbb{Z}_+}\,  x\cdot ||x\alpha_1  || \cdot
||x\alpha_2 -\eta ||
>0.
$$
}

Quite similar result was obtained recently by U. Shapira
\cite{SHA} by means of dynamical systems.
We would like to note here that two papers by E. Lindenstrauss and U. Shapira
\cite{u1,u2} related to the topic appeared very
 recently.

Proof of Corollary 1.

 As $\alpha_1,
\alpha_2$ are linearly dependent, we have integers $A_1, A_2 ,B$,
not all zero, such that
$$
A_1\alpha_1+A_2 \alpha_2 +B = 0.
$$

From Khintchine's Theorem A we can deduce that there exists {\it
uncountably many} $\eta $ satisfying the conclusion of the
theorem. (From \cite{Ts} we know that the corresponding set is a
winning set and hence is uncountable and dense). So we may find
$\eta_1, \eta_2$ satisfying
\begin{equation}\label{ett1}
\inf_{x\in \mathbb{Z}_+}\,  x\cdot ||\alpha_ix - \eta_i|| \ge
\delta,\,\,\,\, i = 1,2.
\end{equation}
and
$$
 ||A_1\eta_1+ A_2\eta_2 ||\ge
\delta
$$
with some positive $\delta$. (For the statement (ii) one can take
$\eta_1 = 0, \eta_2 = \eta$.) Then
\begin{equation}\label{ett2}
 \delta \le ||A_1\eta_1+ A_2\eta_2 ||= ||A_1 (\alpha_1 x - \eta_1) + A_2 (\alpha_2 x - \eta_2) ||
 \le  A\cdot \max_{i=1,2} ||\alpha_ix - \eta_i||, \,\,\,\, A = \max_{i=1,2} |A_i|
 .
\end{equation}
Take a positive integer $x$. From (\ref{ett2}) we see that one of the
quantities $||\alpha_ix - \eta_i|| i = 1,2 $ is not less than
$\delta / A$. To the other quantity we may apply lower bound from
 (\ref{ett2}). This gives
$$
x\cdot ||x\alpha_1  || \cdot ||x\alpha_2 -\eta ||\ge \delta^2/A.
$$
Corollary 1 is proved.

For $\alpha =(\alpha_1, \alpha_2) \in \mathbb{R}^2$ we define a
function
$$
\Psi_\alpha (t) = \min_{(x_1,x_2)\in \mathbb{Z}^2\setminus
\{(0,0)\}, \max |x_i| \le t} \,\, ||\alpha_1x_1+\alpha_2 x_2||.
$$

Now we formulate another two theorems from Khintchine's paper
\cite{k26}.

{\bf Theorem B.} {\it Given a function $\varphi (t)$ decreasing to
zero there exist $\alpha_1, \alpha_2$ linearly independent over
$\mathbb Z$ together with $1$ such that for all $t$ large enough
$$\Psi_\alpha (t) \le \varphi (t).
$$
 }

 {\bf Theorem C.} {\it Given a function $\psi (t)$ increasing to
infinity  there exist  reals $\alpha_1, \alpha_2$ linearly
independent over $\mathbb Z$ together with $1$ and reals
$\eta_1,\eta_2$ such that
$$
\inf _{x\in \mathbb{Z}_+} \psi (x) \cdot \max _{i=1,2} ||\alpha_i
x - \eta_i||
> 0.
$$}

In fact A. Khintchine deduces Theorem C from  Theorem B. In the
fundamental paper  \cite{k26} A. Khintchine  states also two
additional general results. One of them is as follows.

 {\bf Theorem D.} {\it Given a
tuple of real numbers $(\eta_1, \eta_2)$ and given a
  function $\psi (t)$ increasing to
infinity  there exist  reals $\alpha_1, \alpha_2$ linearly
independent over $\mathbb Z$ together with $1$  such that
$$
\inf _{x\in \mathbb{Z}_+} \psi (x) \cdot \max _{i=1,2} ||\alpha_i
x - \eta_i||
> 0.
$$}

From another hand, by a result of J. Tseng \cite{Ts}, we know for any real
$\alpha$   the set
$$
{\cal B} = \{ \eta:\,\, \inf_{x\in \mathbb{Z}_+} x\cdot ||\alpha x
- \eta||>0\}
$$
is an $1/2$-winning set in $\mathbb{R}$. It follows that  the sets
$$
{\cal B} _1 = \{ (\eta_1,\eta_2):\,\,\,  \eta_1 \in {\cal B}, \,
\eta_2 \in \mathbb{R}\},\,\,\,\,\, \,\,\, {\cal B} _2 = \{
(\eta_1,\eta_2):\,\,\, \eta_1 \in \mathbb{R},\,  \eta_2 \in {\cal
B} \ \}
$$
are $1/2$-winning sets in $\mathbb{R}^2$.

 In the paper \cite{MMJ}
N. Moshchevitin proved a general result. The theorem below is a
particular case of this result.

 {\bf Theorem E.} {\it
    Suppose that $\psi
(t) $ is a function increasing to infinity  as $t\to+\infty$.
Suppose that for any $w\ge 1$ we have the inequality
\begin{equation}\label{addi}
\sup_{x\ge 1}\frac{\psi (wx)}{\psi (x)} <+\infty.
\end{equation}
Let $\rho (t)$ be the function inverse to the function
 $t\mapsto 1/\psi (t)$
Let $\alpha=(\alpha_1, \alpha_2)\in \mathbb{R}^2$ be such that
$$
\Psi_\alpha (t) \le \rho (t).
$$
 Then the set
$$
{\cal B}^{[\psi]} = \{ (\eta_1,\eta_2):\,\,\, \inf_{x\in
\mathbb{Z}_+} \psi (x)\cdot \max _{i=1,2} ||\alpha_i x - \eta_i||
> 0\}
 $$
is an $1/2$-winning set in $\mathbb{R}^2$.

 }

From the theory of winning sets (see \cite{WMS}) we know that a
countable intersection of $\alpha$-winning set is an
$\alpha$-winning set also. In particular the set
$$
{\cal B}^{[\psi]} \cap {\cal B}_1\cap {\cal B}_2
$$
is an $1/2$-winning set in $\mathbb{R}^2$.
 Moreover every
$\alpha$-winning set has full Hausdorff dimension and hence is not
empty. Thus  we deduce the following result.

{\bf Theorem 3.}\,\, {\it
 Suppose that $\psi
(t) $ is a function increasing to infinity  as $t\to+\infty$.
Suppose that (\ref{addi}) is valid. Then there exist real numbers
$\alpha_1, \alpha_2$ linearly independent over $\mathbb{Z}$
together with $1$ and real numbers $\eta_1,\eta_2$ such that
$$
\inf_{x\in \mathbb{Z}_+}  x \psi (x)\cdot  ||\alpha_1 x - \eta_1||
\cdot  ||\alpha_2 x - \eta_2|| >0.
$$

}

A proof immediately follows from the fact that ${\cal B}^{[\psi]}
\cap {\cal B}_1\cap {\cal B}_2\neq \varnothing$. Let
$(\alpha_1,\alpha_2) $  be the tuple from  Theorem C applied to
$\varphi (t) = \rho (t)$. Take  $(\eta_1,\eta_2) \in {\cal
B}^{[\psi]} \cap {\cal B}_1\cap {\cal B}_2$. Take positive integer
$x$. One of the values $||\alpha_i x - \eta_i||$ should be greater
than $\varepsilon/\psi (x)$ where $\varepsilon$ depends on
$\alpha_1,\alpha_2,\eta_1,\eta_2$ only. Then the other one is
greater than $\varepsilon '/x$ where $\varepsilon'$ depends on
$\alpha_1,\alpha_2,\eta_1,\eta_2$ only. Theorem 3 is proved.

Theorem 3 may be compared with the main result from the paper
\cite{SHA}. It does not answer the following question, 
already posed in \cite{Bu11}.

{\bf Problem.} {\it
Let $\alpha$ and $\beta$ be real numbers with $1, \alpha, \beta$
being linearly independent over the rationals.
Let $\alpha_0$, $\beta_0$ and $\gamma$ be real numbers.
To prove or to disprove that
$$
\inf_{q \not= 0} \, |q| \cdot \Vert q \alpha - \alpha_0 \Vert \cdot 
\Vert q \beta - \beta_0 \Vert = 0  
$$
and/or that
$$
\inf_{(x, y) \not= (0, 0)} \, \Vert x \alpha + y \beta - \gamma \Vert 
\cdot \max\{\vert x \vert,1\}\cdot\max\{\vert y \vert,1\} = 0. 
$$

}

The following two theorems by U. Shapira from the paper \cite{SHA}
worth noting in the context of this problem.

 {\bf Theorem F.} {\it Almost all (in the sense of Lebesgue measure) pairs
$(\alpha_1,\alpha_2) \in \mathbb{R}^2$ satisfy the following property: for every pair
$(\eta_1,\eta_2) \in \mathbb{R}^2$ one has
$$
\liminf_{q\to \infty} q\,\,||q\alpha_1-\eta_1||\,\,
||q\alpha_1 - \eta_2 || = 0.
$$}

{\bf Theorem G.} {\it 
The conclusion  of Theorem F is true for numbers $\alpha_1, \alpha_2$ which form together with $1$
a basis of a totally real algebraic field of degree $3$.}

Also we would like to refer to one more Khintchine's result (see \cite{k26}, Hilfssatz 4)

 {\bf Theorem H.} {\it
Given $c\in  (0,1)$ there exists $\Gamma>0$ with the following property.
For any $\alpha \in \mathbb{R}$ there exists $\beta \in \mathbb{R}$ such that
$$
\max (cx |\alpha x - y|,\,\, \Gamma |\beta x - z|) \ge 1,
$$
where maximum is taken over integers $x>0,y,z,\,\, (x,y) =1$.
In other words if
$$
 |\alpha x - y| \le \frac{1}{cx}, \,\,\, (x,y) = 1
$$
then
$$
 ||\beta x|| \ge \frac{1}{\Gamma}.
$$

}

At the end of this section we want to refer to wondeful 
recent result by D. Badziahin,  A. Pollington and S. Velani
from the paper \cite{BADZ}. 
In this paper they solve famous W.M. Schmidt's conjecture \cite{schBAD}.

 {\bf Theorem I.} {\it
Let $u,v\ge 0, u+v =1$.
Suppose that
\begin{equation}\label{BPV}
\inf_{x\in \mathbb{Z}_+} 
x^{\frac{1}{u}} ||\alpha x||> 0.
\end{equation}
Then the set
$$
B_{u}(\alpha) =
\{
\beta \in \mathbb{R}:\,\,\,
\inf_{x\in \mathbb{Z}_+}
\, \max (x^u||\alpha x||,\, x^v||\beta x||) \,\, >0
\}
$$
has full Hausdorff dimension.

}

Here we should note that the main result from  \cite{BADZ}
shows for a given $\alpha$ 
under the condition (\ref{BPV})
that intersections of  sets of the form 
$B_{u}(\alpha) $ for a finite collection 
of different values of $u$  has full Hausdorff dimension.
An explicit version of the original proof invented by
D. Badziahin, A. Pollington and S. Velani
  was given in \cite{MoshBBB}, in the simplest case $u=1/2$.

Recently D. Badziahin \cite{BADZ1} proved the following result.

{\bf Theorem J.} {\it
The set
$$
\{(\alpha,\beta )\in \mathbb{R}^2:\,\,\, \inf_{x\in \mathbb{Z}, x \ge 3}
x\log x \,\log\log x  \,||\alpha x|| \,
||\beta x|| >0\}
$$
has Hausdorff dimension equal to 2.

Moreover 
if $\alpha $ is a badly approximable number then the set
$$
\{\beta \in \mathbb{R}:\,\,\, \inf_{x\in \mathbb{Z}, x \ge 3}
x\log x \,\log\log x  \,||\alpha x|| \,
||\beta x|| >0\}
$$
has Hausdorff dimension equal to 1.
}

We think that the method from \cite{BADZ,BADZ1} cannot be generalized for inhomogeneous setting.

 {\bf 4. Examples to Theorem 1.}
Here we give several special choices of   parameters in Theorem 1
and deduce several corollaries.

{\bf Example 1.} Put
$$
\omega_1(t) = \omega _2(t) = \gamma t
$$
with some positive $\gamma >1$. Then
$$
\omega_1^* (t) = \frac{t}{\gamma}
$$
and we may take in  (\ref{amega})
$$
\Omega (y,z) = \sqrt{\frac{1}{\gamma}\, \frac{z}{y}} .
$$
Put
$$
\psi_1(t) = \psi_2 (t) = t, \,\,\,\,\, \phi (t) = t\cdot  {\ln
^2t} .
$$
Then
$$
\psi_2^* (t) = t.
$$
So
\begin{equation}\label{ddl}
\delta^{[1]}_\varepsilon (\mu,\nu)= 2 \cdot \varepsilon \cdot
\frac{2^{\mu - \nu}}{\nu^2}
\end{equation}
and
$$
 S_{A,\varepsilon}^{[1]}(X)=\,\,
2  \cdot \varepsilon\cdot
 \,\sum_{X\le \nu <A(X+1)}\, \sum_{1\le \mu \le
\nu + \log_2\gamma + 2} \frac{2^{\mu - \nu}}{\nu^2}
 \max\left(\sqrt{\frac{1}{\gamma}\,2^{\nu - \mu +2}} , 2^{\nu  -\mu},1 \right)
 \le
 $$
 $$\le
4 \cdot \varepsilon\cdot \sum_{X\le \nu <A(X+1)}\,\left(
\sum_{1\le \mu \le \nu } \frac{1}{\nu^2} + \sum_{\nu+1\le \mu \le
\nu +2 +\log_2\gamma} \frac{2^{\mu - \nu}}{\nu^2}
 \right)\le 8  \cdot \varepsilon\cdot \sum_{X\le \nu
 <A(X+1)}\,\left( \frac{1}{\nu } + \frac{4\gamma
 }{\nu^2}\right)\le
$$
$$
\le 16 \varepsilon \ln (2A)
$$
for $X_0$ large enough ($X_0 \ge \gamma   /\varepsilon$). Put
$A=4$. Then the condition (\ref{cond0}) is satisfied provided
$\frac{X_0}{\ln^2 X_0} \ge \frac{1}{\varepsilon}$. Thus we obtain
the following results.

{\bf Corollary 1.1.} {\it Let $\eta_x, x =1,2,3,..$ be a  sequence
of reals. Given positive $\varepsilon \le 2^{-14}$ and a badly
approximable real $\alpha$ such that
$$
||\alpha x || \ge \frac{1}{\gamma x}\,\,\, \forall x \in
\mathbb{Z}_+, \,\,\, \gamma >1,$$ there exist $X_0 = X_0
(\varepsilon ,\gamma )$ and a real $\beta$ such that
$$
\inf_{x\ge X_0}\,  x \ln^2x \cdot ||x\alpha|| \cdot ||x\beta
-\eta_x||\, \ge \varepsilon.
$$

}

{\bf Corollary 1.2.} {\it Let $\eta_x, x =1,2,3,..$ be a  sequence
of reals. Given positive $\varepsilon \le 2^{-14}$ and  real
$\alpha, \eta$ such that simultaneously
$$
||\alpha x || \ge \frac{1}{\gamma x}\,\,\, \forall x \in
\mathbb{Z}_+, \,\,\, \gamma >1$$ and
$$
||\alpha x - \eta|| \ge \frac{1}{\gamma x}\,\,\, \forall x \in
\mathbb{Z}_+, \,\,\, \gamma >1,$$ there exist $X_0 = X_0
(\varepsilon ,\gamma )$ and a real $\beta$ such that
$$
\inf_{x\ge X_0}\,  x \ln^2x \cdot ||x\alpha-\eta|| \cdot ||x\beta
-\eta_x||\, \ge \varepsilon.
$$

}

From Khintchine's Theorem A we deduce the following result.

{\bf Corollary 1.3.} {\it Let $\eta_x, x =1,2,3,..$ be a  sequence
of reals. Given positive $\varepsilon \le 2^{-14}$ and a   real
$\alpha$ such that
$$
||\alpha x || \ge \frac{1}{\gamma x}\,\,\, \forall x \in
\mathbb{Z}_+, \,\,\, \gamma >1,$$ there exist $X_0 = X_0
(\varepsilon ,\gamma )$ and  real $ \eta , \beta$ such that
$$
\inf_{x\ge X_0}\,  x \ln^2x \cdot ||x\alpha-\eta || \cdot ||x\beta
-\eta_x||\, \ge \varepsilon.
$$

}

{\bf Example 2.} Put
$$
\omega_1(t) = \omega _2(t) =  t \ln t
$$
 Then
$$
\omega_1^* (t) \asymp \frac{t}{ \ln t}
$$
and we may take in  (\ref{amega})
$$
\Omega (y,z) = c\sqrt{ \, \frac{z \ln z}{y}}
$$
with  small positive $c$.

Put
$$
\psi_1(t) = \psi_2 (t) = t, \,\,\,\,\, \phi (t) = t\cdot  {\ln
^2t} .
$$
Then
$$
\psi_2^* (t) = t ,
$$
and again $ \delta^{[1]}_\varepsilon (\mu,\nu)$ satisfies
(\ref{ddl}). Now
$$
 S_{A,\varepsilon}^{[1]}(X) \ll
\varepsilon\cdot
 \,\sum_{X\le \nu <A(X+1)}\, \sum_{1\le \mu \le
\nu + \log_2(\nu +1) + 2} \frac{2^{\mu - \nu}}{\nu^2}
 \max\left(\sqrt{2^{ {\nu - \mu} }\nu} , 2^{\nu  -\mu},1 \right)
 \ll
 $$
 $$
 \ll
\varepsilon\cdot
 \,\sum_{X\le \nu <A(X+1)}\, \sum_{1\le \mu \le
\nu + \log_2(\nu +1) + 2}  \frac{2^{\mu - \nu}}{\nu^2}
 \max\left(\sqrt{  2^{\nu - \mu }\nu} , 2^{\nu  -\mu} \right)
 \ll
 $$
 $$
 \ll
  \varepsilon\cdot \sum_{X\le \nu <A(X+1)}\,\left(
\sum_{1\le \mu \le \nu } \frac{1}{\nu^2} + \sum_{\nu - \log_2
(\nu+1)\le \mu \le \nu  +\log_2(\nu+1) + 2} \frac{2^{\frac{\mu -
\nu}{2}}}{\nu^{3/2}}\right)\ll \varepsilon \ln 2A,
$$
for $X_0$ large enough. Put $A=4$. Then for $X_0$ large enough the
inequality  (\ref{cond0}) is valid. Thus we obtain the following
results.

{\bf Corollary 2.1.} {\it There exists an absolute positive
constant $\varepsilon_0 $ with the following property. Let
$\eta_x, x =1,2,3,..$ be a  sequence of reals. Given positive
$\varepsilon \le \varepsilon_0$ and a  real $\alpha$ such that for
all $x\ge X_1$ one has
$$
||\alpha x || \ge \frac{1}{ x\ln x}\ ,  $$ there exist $X_0 = X_0
(\varepsilon ,X_1 )$ and a real $\beta$ such that
$$
\inf_{x\ge X_0}\,  x \ln^2x \cdot ||x\alpha|| \cdot ||x\beta
-\eta_x||\, \ge \varepsilon.
$$

}

Corollary 2.1 is a more general statement than Corollary 1.1.

{\bf Corollary 2.2.} {\it Let $\eta_x, x =1,2,3,..$ be a  sequence
of reals. Given for positive $\varepsilon $ small enough and  real
$\alpha, \eta$ such that for all $ x\ge X_1$ simultaneously
$$
||\alpha x || \ge \frac{1}{x\ln x} , $$ and
$$
||\alpha x - \eta|| \ge \frac{1}{ x\ln x} , \,$$ there exist $X_0
= X_0 (\varepsilon ,X_1 )$ and a real $\beta$ such that
$$
\inf_{x\ge X_0}\,  x \ln^2x \cdot ||x\alpha-\eta|| \cdot ||x\beta
-\eta_x||\, \ge \varepsilon.
$$

}

{\bf Example 3.} Put
$$
\omega_1(t) =    t \ln^2 t,\,\,\,\,\ \omega _2(t) =\gamma t,
\gamma
>1.
$$
 Then
$$
\omega_1^* (t) \asymp \frac{t}{ \ln^2 t}
$$
and we may take in  (\ref{amega})
$$
\Omega (y,z) = c\sqrt{ \, \frac{z  }{y}}\ln z
$$
with  small positive $c$.

Put
$$
\psi_1(t) = \psi_2 (t) = t, \,\,\,\,\, \phi (t) = t\cdot  {\ln
^2t} .
$$
Then
$$
\psi_2^* (t) = t ,
$$
and again $ \delta^{[1]}_\varepsilon (\mu,\nu)$ satisfies
(\ref{ddl}). So
$$
 S_{A,\varepsilon}^{[1]}(X) \ll
\varepsilon\cdot
 \,\sum_{X\le \nu <A(X+1)}\, \sum_{1\le \mu \le
\nu + \log_2\gamma + 2} \frac{2^{\mu - \nu}}{\nu^2}
 \max\left(\sqrt{2^{ {\nu - \mu} } } \, \cdot \, \nu, 2^{\nu  -\mu}  \right)
 \ll
 $$
 $$
 \ll
  \varepsilon\cdot \sum_{X\le \nu <A(X+1)}\,\left(
\sum_{1\le \mu \le \nu } \frac{1}{\nu^2} + \sum_{\nu - 2\log_2
(\nu+1)\le \mu \le \nu  + \log_2\gamma + 2} \frac{2^{\frac{\mu -
\nu}{2}}}{\nu }\right) \ll
$$
$$
\ll
 \varepsilon\cdot \sum_{X\le \nu <A(X+1)}\,
\sum_{1\le \mu \le \nu } \frac{1+\sqrt{\gamma}}{\nu} \ll
\varepsilon (1+\sqrt{\gamma})\ln 2A,
$$
for $X_0$ large enough. Again with  $A=4$ for $X_0$ large enough
the inequality  (\ref{cond0}) is valid. Thus we obtain the
following results. This result  is a more general statement than
Corollary 1.2.

{\bf Corollary 3.1.} {\it Let $\eta_x, x =1,2,3,..$ be a sequence
of reals. Let $\gamma >1$. Suppose that the product $\varepsilon
\sqrt{\gamma} $ is small enough. Suppose that for certain   real
$\alpha, \eta$ and for   $ x\ge X_1$ simultaneously one has
$$
||\alpha x || \ge \frac{1}{ x\ln^2 x}$$ and
$$
||\alpha x - \eta|| \ge \frac{1}{\gamma x} .$$ Then there exist
$X_0 = X_0 (\varepsilon, \gamma  ,X_1 )$ and a real $\beta$ such
that
$$
\inf_{x\ge X_0}\,  x \ln^2x \cdot ||x\alpha-\eta|| \cdot ||x\beta
-\eta_x||\, \ge \varepsilon.
$$

}

Now from Khintchine's Theorem A we deduce a result which is more
general that Corollary 1.3.

{\bf Corollary 3.2.} {\it Let $\eta_x, x =1,2,3,..$ be a  sequence
of reals. Given positive $\varepsilon  $ small enough and a   real
$\alpha$ such that
$$
||\alpha x || \ge \frac{1}{ x \ln^2 x } ,$$ there exist $X_0 = X_0
(\varepsilon   )$ and  real $ \eta , \beta$ such that
$$
\inf_{x\ge X_0}\,  x \ln^2x \cdot ||x\alpha-\eta || \cdot ||x\beta
-\eta_x||\, \ge \varepsilon.
$$

}

{\bf Example 4.} Put
$$
\omega_1(t) = \omega _2(t) = \gamma t
$$
with some positive $\gamma >1$. Then as in Example 1 we have
$$
\omega_1^* (t) = \frac{t}{\gamma},\,\,\,\,\,\ \Omega (y,z) =
\sqrt{\frac{1}{\gamma}\, \frac{z}{y}} .
$$
Suppose that $0\le a < 1$. Put
$$
\psi_1(t) =  t \cdot (\log_2 1/ t)^a,\,\,\,\,\psi_2 (t) = t,
\,\,\,\,\, \phi (t) = t\cdot {\log_2 ^{2-a} t} .
$$
Then
$$
\psi_2^* (t) = t
$$
and
$$
\delta^{[1]}_\varepsilon (\mu,\nu) = 2 \cdot \varepsilon \cdot
\frac{2^{\mu - \nu}}{\nu^{2-a} (\mu+1)^a}.
$$
Now
$$
 S_{A,\varepsilon}^{[1]}(X) \ll
4\varepsilon\cdot
 \,\sum_{X\le \nu <A(X+1)}\, \sum_{1\le \mu \le
\nu + \log_2\gamma + 2} \frac{2^{\mu - \nu}}{\nu^{2-a} (\mu+1)^a}
 \max\left(\sqrt{2^{ {\nu - \mu} }/\gamma } \, , 2^{\nu  -\mu}, 1  \right)
 \le
 $$
 $$
 \le
4 \cdot \varepsilon\cdot \sum_{X\le \nu <A(X+1)}\,\left(
\sum_{1\le \mu \le \nu } \frac{1}{\nu^{2-a} (\mu+1)^a} +
\sum_{\nu+1\le \mu \le \nu +2 +\log_2\gamma} \frac{2^{\mu -
\nu}}{\nu^{2-a} (\mu+1)^a}
 \right)\le
 $$
 $$
 \le
 \frac{8\varepsilon}{1-a}
\sum_{X\le \nu <A(X+1)}\left(\frac{1}{\nu} + \frac{4\gamma
}{\nu^2}
 \right)
 \le \frac{32\varepsilon \ln (2A)}{1-a}
$$
for  $X_0 \ge \gamma   /\varepsilon$. Put again $A=4$. Then the
condition (\ref{cond0}) is satisfied provided $\frac{X_0}{\ln^2
X_0} \ge \frac{1}{\varepsilon}$. Thus we obtain the following
results (compare with Theorem 3 from \cite{MB}).

{\bf Corollary 4.1.} {\it Let $\eta_x, x =1,2,3,..$ be a  sequence
of reals. Suppose that $0\le a <1$. Given positive $\varepsilon
\le\frac{1}{2^{20}(1-a)}$ and a badly approximable real $\alpha$
such that
$$
||\alpha x || \ge \frac{1}{\gamma x}\,\,\, \forall x \in
\mathbb{Z}_+, \,\,\, \gamma >1,$$ there exist $X_0 = X_0
(\varepsilon ,\gamma )$ and a real $\beta$ such that
$$
\inf_{x\ge X_0}\,  x (\log_2  x )^{2-a} \cdot (\log_2 1/
||x\alpha||)^a \cdot
  ||x\alpha|| \cdot ||x\beta -\eta_x||\, \ge \varepsilon.
$$

}

{\bf Corollary 4.2.} {\it Let $\eta_x, x =1,2,3,..$ be a  sequence
of reals. Given positive $\varepsilon 
\le\frac{1}{2^{20}(1-a)}$ and  real $\alpha, \eta$ such that
simultaneously
$$
||\alpha x || \ge \frac{1}{\gamma x},\,\,\,\,\, ||\alpha x -
\eta|| \ge \frac{1}{\gamma x}\,\,\, \forall x \in \mathbb{Z}_+,
\,\,\, \gamma >1,$$ there exist $X_0 = X_0 (\varepsilon ,\gamma )$
and a real $\beta$ such that
$$
\inf_{x\ge X_0}\,  x (\log_2  x )^{2-a} \cdot (\log_2 1/
||x\alpha-\eta||)^a \cdot
  ||x\alpha-\eta || \cdot ||x\beta -\eta_x||\, \ge \varepsilon.
$$

}

 {\bf Corollary 4.3.} {\it Let $\eta_x, x =1,2,3,..$ be a  sequence
of reals. Given positive $\varepsilon \le\frac{1}{2^{20}(1-a)}$
and a   real $\alpha$ such that
$$
||\alpha x || \ge \frac{1}{\gamma x}\,\,\, \forall x \in
\mathbb{Z}_+, \,\,\, \gamma >1,$$ there exist $X_0 = X_0
(\varepsilon ,\gamma )$ and  real $ \eta , \beta$ such that
$$
\inf_{x\ge X_0}\,  x (\log_2  x )^{2-a} \cdot (\log_2 1/
||x\alpha-\eta ||)^a \cdot
  ||x\alpha-\eta || \cdot ||x\beta -\eta_x||\, \ge \varepsilon.
$$

}

Of course one can deduce other corollaries of a similar type from
Theorem 1. For example one may deduce statements which are more
general than Corollaries 4.1 - 4.3 in the same manner as it was
done in Examples 2,3.

 {\bf 5. Examples to Theorem 2.}
Here we consider some corollaries related to   special choices of
parameters in Theorem 2.

{\bf Example 5.}

Let $u, v > 0,\,\,\, u+v =1$.
Put
$$
\omega_1(t)  =  \frac{\gamma  t^{\frac{1}{u}}}{(\ln t)^u},\,\,\,\gamma
>1.
$$
 Then
  we may take in  (\ref{amega})
$$
\Omega (y,z) = c \left(\frac{z  }{y^u (\ln z)^{u^2} }\right)^{\frac{1}{1+u}} $$
with  small positive $c$ (we take into account that  $x\ll z^{1/2} \ln z$) .

Put
$$
\psi_1(t) = \psi_2 (t) = t, 
\,\,\,\,\,
\phi_1(t) = (t\log_2 t)^u,\,\,\, \phi_2 (t) = (t\log_2 t)^v. 
$$
Then
$$
\psi_1^* = 
\psi_2^* (t) = t ,
$$
and  
$$ \delta^{[2]}_\varepsilon (\nu)
=\frac{\varepsilon}{(\nu 2^\nu)^v},\,\,\,\,\, r_\varepsilon (\nu) = \frac{\varepsilon}{(\nu 2^\nu)^u}.
$$ 
 So
$$
 S_{A,\varepsilon}^{[2]}(X) \ll \varepsilon
\cdot
 \,\sum_{X\le \nu <A(X+1)}\,  \frac{1}{(\nu 2^\nu)^v} \cdot \frac{2^{(1-u)\nu} }{\nu^u}
 \ll
  \varepsilon\cdot \sum_{X\le \nu <A(X+1)}\,
\sum_{1\le \mu \le \nu } \frac{1 }{\nu} \ll
\varepsilon \ln 2A,
$$
 So we get

{\bf Corollary 5.1.} {\it
Suppose that  $u, v > 0, \,\, u+v = 1$. 
 Let $\eta_x, x =1,2,3,..$ be a sequence
of reals.
Let $\eta$ be an arbitrary real number.
 Let $\gamma >0$. Suppose that  $\varepsilon
$ is small enough. Suppose that for certain   real
$\alpha$ and for   $ x\ge X_1$  one has
$$
||\alpha x || \ge \frac{\gamma \,(\ln x)^u}{ x^{1/u}}.$$   Then there exist
$X_0 = X_0 (\varepsilon, \gamma  ,X_1 )$ and a real $\beta$ such
that
$$
\inf_{x\ge X_0}\,  \max  ((x \ln x)^u \cdot ||x\alpha-\eta||, \,\,(x\ln x)^v \cdot ||x\beta
-\eta_x|| )\, \ge \varepsilon.
$$

}

{\bf Example 6.}

Put
$$
\omega_1(t)  =  {\gamma t}{\ln t},\,\,\,\gamma
>1.
$$
 Then
  we may take in  (\ref{amega})
$$
\Omega (y,z) = c\sqrt{ \, \frac{z \ln z}{y}}
$$
with  small positive $c$.

Put
$$
\psi_1(t) = \psi_2 (t) = t, 
\,\,\,\,\,
\phi_1(t) = \Delta t,\,\,\, \phi_2 (t) = (\log_2 t)^{3/2},\,\,\,\,
\Delta >0.
$$
Then
$$
\psi_1^* = 
\psi_2^* (t) = t ,
$$
and  
$$ \delta^{[2]}_\varepsilon (\nu)
=\frac{\varepsilon}{\nu^{3/2} },\,\,\,\,\, r_\varepsilon (\nu) = \frac{\varepsilon}{\Delta 2^\nu}.
$$ 
 So
$$
 S_{A,\varepsilon}^{[2]}(X) \ll \varepsilon
\cdot
 \,\sum_{X\le \nu <A(X+1)}\,  \frac{\Omega (1/2^{\nu+1})}{\nu^{3/2}}\ll
 \frac{ \varepsilon^{3/2}}{\Delta^{1/2}}\cdot \sum_{X\le \nu <A(X+1)}\,
\sum_{1\le \mu \le \nu } \frac{1 }{\nu} \ll
\frac{ \varepsilon^{3/2}}{\Delta^{1/2}} \ln A,
$$
 So we get

{\bf Corollary 6.1.} {\it
 Let $\eta_x, x =1,2,3,..$ be a sequence
of reals.
Let $\eta$ be an arbitrary real number.
 Let $\Delta >1$. Suppose that  $2^{20}\varepsilon^3\le \Delta
$ . Suppose that for certain   real
$\alpha$ and for  all positive integers $ x$  one has
$$
||\alpha x || \ge \frac{\gamma}{ x\ln x}.$$   Then there exist
$X_0 = X_0 (\gamma  )$ and a real $\beta$ such
that
$$
\inf_{x\ge X_0}\,  \max  (\Delta  x \cdot ||x\alpha-\eta||,\,\, (\ln x)^{3/2} \cdot ||x\beta
-\eta_x|| )\, \ge \varepsilon.
$$
In other words for  this $\beta$ if
$$
||x\alpha-\eta||\le \frac{\varepsilon}{\Delta x}
$$
then
$$
||x\beta
-\eta_x||\ge \frac{\varepsilon}{(\ln x)^{3/2}}.
$$

}

{\bf 6. Sets of integers.}

Consider sets
$$
A_{\nu, \mu} = \{ x\in \mathbb{Z}_+:\,\, 2^\nu \le
x<2^{\nu+1},\,\,2^{-\mu -1}< ||\alpha x-\eta||\le 2^{-\mu}\},
$$
$$
A_{\nu } (t) = \{ x\in \mathbb{Z}_+:\,\, 2^\nu \le
x<2^{\nu+1},\,\, ||\alpha x-\eta||\le t\},
$$
Now we deduce an upper bound for the cardinality of the set
$A_{\nu, \mu}$.

{\bf Lemma 1.}\,\,\,{\it  Under the condition (\ref{bad})
one has
$$
{\rm card}\,A_{\nu, \mu} \le  2^3 \max\left(\Omega
(2^{\mu-1}, 2^{\nu +1}), 2^{\nu  -\mu} ,1\right).
$$
}

Proof.
 For $ a\in A_{\nu,\mu}$ define integer $y$ from
the condition
$$
||x\alpha -\eta||= |x\alpha -\eta-y| .$$

{\bf  Case 1$^0$.} All integer points $z= (x,y), x\in A_{\nu,
\mu}$ form a convex polygon $\Pi$ of positive measure ${{\rm mes}
\, \Pi}>0 $. Then
\begin{equation}\label{01}
{\rm card}\,  A_{\nu, \mu} \le 6 \, {{\rm mes} \, \Pi} \le 6\cdot
2^{\nu +1 -\mu }< 2^{\nu  -\mu+3}.
\end{equation}

{\bf  Case 2$^0$.} All integer points $z= (x,y), x\in A_{\nu,
\mu}$ lie on the same line. Then all these points are of the form
$$
z_0 + l z_1,\,\,\, z_j = (x_j,y_j),\,\,\, 0\le l \le L.
$$
Now we see that
$$
|\alpha Lx_1 - Ly_1| \le 2^{-\mu+1}
$$
and
$$
|\alpha x_1 - y_1| \le 2^{-\mu+1}L^{-1}.
$$
From (\ref{bad}) we have
$$
\omega_1 (x_1) \ge 2^{\mu-1} L.
$$
So
$$
x_1\ge \omega_1^* (2^{\mu-1} L)
$$
and
$$
Lx_1\le 2^{\nu+1}.
$$
We conclude that
$$
L\le 2^{\nu+1},\,\,\,\, L\cdot \omega_1^* (2^{\mu-1} L) \le
2^{\nu+1 }.
$$
So by (\ref{amega00}) we have
\begin{equation}\label{02}
{\rm card}\,A_{\nu, \mu} \le L+1 \le \Omega (2^{\mu-1}, 2^{\nu
+1})+1.
\end{equation}
We take together (\ref{01},\ref{02}) to obtain
$$ 
{\rm card}\,A_{\nu, \mu} \le \max\left(\Omega (2^{\mu-1},  2^{\nu +1}), 2^{\nu  -\mu},1 \right).
$$
Lemma is proved.

The next 
lemma deals with the cardinality of $A_\nu (t)$.

{\bf Lemma 2.}\,\,\,{\it  Under the condition (\ref{bad})
one has
$$
{\rm card}\,A_{\nu} (t) \le  2^2 \max\left(\Omega
(1/2t, 2^{\nu +1}), 2^{\nu }t ,1\right).
$$
}

Proof.

The proof is quite similar to the proof of Lemma 1.
We should consider two similar cases 1$^0$ and 2$^0$.
In the {\bf Case 1$^0$} we deduce the bound
$${\rm card}\,A_{\nu} (t) \le 2^{\nu+2}t.$$
In the {\bf Case 2$^0$} we see that 
$$
L\le 2^{\nu+1},\,\,\,\, L\cdot \omega_1^* \left(\frac{ L}{2t}\right) \le
2^{\nu+1 }.
$$
By (\ref{amega00}) we have
$$
{\rm card}\,A_{\nu, \mu} \le L+1 \le \Omega (1/2t, 2^{\nu
+1})+1.
$$
Lemma 2 follows.

{\bf 7. Lemmas about fractional parts.}

 Put
\begin{equation}\label{delta}
\sigma^{[1]}_\varepsilon (x) = \sigma^{[1]}_{\varepsilon,
\alpha,\gamma} (x)= \psi_2^*\left( \frac{ \varepsilon}{\phi (x)
\psi_1 (||x\alpha-\gamma||)}
 \right),
\end{equation}
\begin{equation}\label{delta2}
\sigma^{[2]}_\varepsilon (x) = \sigma^{[2]}_{\varepsilon,
\alpha,\gamma} (x)= \psi_2^*\left( \frac{ \varepsilon}{\phi_2 (x)
}
 \right).
\end{equation}

Then from the definitions (\ref{delta}) of
$\sigma^{[1]}_\varepsilon (x) $ and $\delta^{[1]}_\varepsilon
(\mu,\nu)$ and monotonicity conditions we see that
\begin{equation}\label{arrow}
x \in A_{\nu, \mu} \,\,\Longrightarrow\,\,
\sigma^{[1]}_\varepsilon (x)\le \delta^{[1]}_\varepsilon
(\mu,\nu).
\end{equation}

Consider  sums
\begin{equation}\label{docond1}
  T^{[1]}_{A,\varepsilon} (Y) = \,\,\,\sum_{Y\le x <Y^A}\, \sigma^{[1]}_\varepsilon (x)
,
\end{equation}
(with $\sigma$ defined in (\ref{delta})) amd
\begin{equation}\label{docond10}
  T^{[2]}_{A,\varepsilon} (Y) = \,\,\,\sum_{Y\le x <Y^A,\,\,
 \phi_1 (x) \psi_1(||\alpha x||) \le \varepsilon
  }\, \sigma^{[2]}_\varepsilon (x)
,
\end{equation}

 {\bf Lemma 3.}

 {\it
Suppose that  (\ref{bad}) and  (\ref{bad1})  are valid. Then under
the condition (\ref{cond})  one has
\begin{equation}\label{cond1}
 \sup_{Y\in \mathbb{Z}_+} \,\,\, T^{[1]}_{A,\varepsilon} (Y)\le \frac{1}{2^6}.
\end{equation}

 }

Proof. Put $X =[\log_2 Y]$.
 We see that
$$
 T_{A,\varepsilon}^{[1]} (Y) \le
 \sum_{X\le \nu < A(X+1)}\,\,
 \sum_{\mu =1}^\infty
\,\, \sum_{ x \in A_{\nu,\mu}} \sigma^{[1]}_\varepsilon (x).
$$
Note that from (\ref{bad1}) it follows that sets $ A_{\nu, \mu}$
are empty for $\mu > \log_2(\omega_2 (2^{\nu+1})) + 1$. So
  from (\ref{arrow}) we have
 \begin{equation}\label{middle}
T_{A,\varepsilon}^{[1]}  (Y) \le
 \sum_{X\le \nu < A(X+1)}
 \sum_{\mu =1}^{[\log_2(\omega_2 (2^{\nu+1}))] + 1}
\,\,\, \delta^{[1]}_\varepsilon (\mu,\nu) \,\times\, {\rm card}\,
A_{\nu, \mu}.
\end{equation}

Now from  (\ref{middle}) and Lemma 1 we have
$$
T_{A,\varepsilon}^{[1]}  (Y) \le
 2^3\, \sum_{X\le \nu < A(X+1)}
 \sum_{\mu =1}^{[\log_2(\omega_2 (2^{\nu+1}))] + 1}
\,\,\, \delta^{[1]}_\varepsilon (\mu,\nu)  \,\times\,
   \max\left(\Omega (2^{\mu-1},
2^{\nu +1}), 2^{\nu  -\mu} ,1 \right).
$$
Lemma 3 follows from (\ref{cond}).

{\bf Lemma 4.}

 {\it
Suppose that  (\ref{bad}) is valid. Then under
the condition (\ref{condmore})  one has
\begin{equation}\label{cond1}
 \sup_{Y\in \mathbb{Z}_+} \,\,\, T^{[2]}_{A,\varepsilon} (Y)\le \frac{1}{2^6}.
\end{equation}

 }

Proof.

The proof is quite similar to those of Lemma 3. Put $X =[\log_2
Y]$. Then 
$$
 T_{A,\varepsilon}^{[2]} (Y) \le
 \sum_{X\le \nu < A(X+1)}\,
\,\, \sum_{ x \in A_{\nu}(r_\varepsilon(\nu ))} \sigma^{[2]}_\varepsilon (x),$$
where $r_\varepsilon (\nu)$ is defined in (\ref{docond3}).
 Now Lemma 4 immediately follows from (\ref{docond2}, \ref{condmore}), Lemma 2 and
the inequality $\sigma^{[2]}_\varepsilon (x) \le  
 \delta^{[2]}_\varepsilon (\nu) $ which is valid for $ x\in A_{\nu}(r_\varepsilon(\nu ))$.

 {\bf 8. Common  PS argument.}
Here we  follow the arguments from the paper \cite{PS} by Y. Peres
and W. Schlag.

Let $j \in \{ 1,2\}$.
 For integers $ 2\le x, 0\le y\le
x$  define
\begin{equation}
E^{[j]} (x,y) = \left[ \frac{y+  \eta_x}{x} -\frac{
\sigma^{[j]}_\varepsilon (x)}{x}, \frac{y+  \eta_x}{x} +\frac{
\sigma^{[j]}_\varepsilon (x)  }{ x} \right],
 \,\,\,
E^{[j]} (x) =\bigcup_{y =0}^{x} E^{[j]} (x,y)\bigcap [0,1].
\label{E}
\end{equation}
 Define
\begin{equation}
l_0 = 0,\,\,\,
 l_x = l_x^{[j]}= [\log_2 (  x /2 \sigma^{[j]}_\varepsilon (x) )  ] ,\,\, x \in \mathbb{N} . \label{L}
\end{equation}
 Each
segment form the union $E_\alpha (x)$ from (\ref{E}) can be
covered by a dyadic interval of the form
$$
\left( \frac{b}{2^{l_x }}, \frac{b+z}{2^{l_x }}\right),\,\,\, z =
1,2 .$$

Let $A^{[j]} (x)$ be the smallest union of all such dyadic
segments which cover the whole set  $E^{[j]} (x)$. Put
$$
(A^{[j]})^c (x) = [0,1] \setminus A^{[j]} (x). $$ Then
$$
(A^{[j]})^c (x)
 = \bigcup_{\nu = 1}^{\tau_x } I_\nu
$$
where closed  segments $I_\nu $ are of the form
\begin{equation}
\left[ \frac{a}{2^{l_x}}, \frac{a+1}{2^{l_x}}\right] ,\,\,\, a\in
\mathbb{Z}. \label{aaa}
\end{equation}

 We take $q_0$ to be a large  positive integer. In order to prove Theorem 1 it is sufficient to show that for all $
q \ge q_0 $ the sets
$$
B^{[1]}_q =\bigcap_{x=q_0}^q (A^{[1]})^c (x)$$ are not empty.
Indeed as the sets $B_q^{[1]}$ are closed and nested we see that
there exists real $\beta$ such that
$$
\beta \in \bigcap_{q\ge q_0} B_q^{[1]} .$$ One can see that the
pair $\alpha, \beta $ satisfies the conclusion of Theorem 1.

Similarly, in order to prove Theorem 2 it is sufficient to show
that for all $ q \ge q_0 $ the sets
$$
B^{[2]}_q =\bigcap_{ x\le q,\,\,\, \phi_1 (x) \psi_1(||\alpha x
||)\le \varepsilon }  (A^{[2]})^c  (x)$$ are not empty.

Under the conditions of Theorems 1 and 2 the following statement
is valid:

{\bf Lemma 5.}\,\,{\it Let $ j \in \{ 1, 2\}$. Suppose that
$\varepsilon $ is small enough. Then for $q_0$ large enough and
for any
$$q_1\ge q_0
,\,\,\, q_2 = q_1^A,\,\,\, q_3 =q_2^A
$$ the following holds. If
\begin{equation}
{\rm mes} B_{q_2}^{[j]} \ge {\rm mes} B_{q_1}^{[j]}/2>0 \label{m1}
\end{equation}
then
\begin{equation}
{\rm mes} B_{q_3}^{[j]} \ge {\rm mes} B_{q_2}^{[j]}/2>0.
\label{m2}
\end{equation}
} Theorems 1, 2  follow from Lemma 5 by induction as the base of
the induction  obviously follows from the arguments of Lemma's
proof.

Proof of  Lemma 5.   First of all we show that for every $j \in
\{1. 2\}$ and  $x\ge q^A$ where $q\ge q_0$ one has
\begin{equation}
 {\rm mes }\left( B_q^{[j]}  \bigcap A^{[j]}   (x) \right) \le 2^4\sigma^{[1]}_\varepsilon (x) \times {\rm mes} B_q^{[j]} .
\label{PP}
\end{equation}
Indeed as from (\ref{L}) and from (\ref{cond0}) in the case $j =
1$ (or from (\ref{cond000}) in the case $ j = 2$) it follows that
$$l_x^{[j]} \le (A-1)\log q ,\, \forall x\le q.$$
We see that
 $B_q^{[j]} $ is a union
$$
B_q^{[j]}  = \bigcup_{\nu = 1}^{T_q } J_\nu
$$
with $J_\nu$ of the form $$\left[ \frac{a}{2^{l}},
\frac{a+1}{2^{l}}\right] ,\,\,\, a\in \mathbb{Z}.
$$
Note that $A^{[j]} (x)$ consists of the segments of the form
(\ref{aaa}) and for $ x\ge q^A > 2^{l+1}
 $ (for $q_0$ large enough)   we see that each $J_\nu$ has at least two rational fractions of the form $\frac{y}{x}, \frac{y+1}{x}$ inside. So
 \begin{equation}
{\rm mes} (J_\nu \cap A^{[j]} (x)) \le 2^4\sigma^{[j]}_\varepsilon
(x)  \times {\rm mes} J_\nu. \label{aaaa}\end{equation} Now
(\ref{PP}) follows from (\ref{aaaa}) by summation over $ 1\le
\nu\le T_q$.

To continue we observe that
$$
B_{q_3}^{[1]} = B_{q_2}^{[1]} \setminus
 \left(\bigcup_{x=q_2+1}^{q_3}  A^{[1]} (x) \right)   ,
$$
and
$$
B_{q_3}^{[2]} = B_{q_2}^{[2]} \setminus
 \left(\bigcup_{q_2+1\le x\le q_3,\, \phi_1 (x) \psi_1(||\alpha x||) \le \varepsilon
    }\,\,\,\,  A^{[2]} (x) \right).
$$
Hence
$$
{\rm mes}  B_{q_3}^{[1]} \ge {\rm mes}   B_{q_2}^{[1]}  -
\sum_{x=q_2+1}^{q_3} {\rm mes} ( B_{q_2}^{[1]}\cap A^{[1]} (x) ).
$$
At the same time
$$
{\rm mes}  B_{q_3}^{[2]} \ge {\rm mes}   B_{q_2}^{[2]}  -
\sum_{q_2+1\le x\le q_3,\, \phi_1 (x) \psi_1(||\alpha x||) \le
\varepsilon
    }\,\,\,\, {\rm mes} ( B_{q_2}^{[2]}\cap A^{[2]} (x) ).
$$

 As
 $$
B_{q_2}^{[j]}\cap A^{[j]} (x)\subseteq B_{q_1}^{[j]}\cap A^{[j]}
(x) $$ we can apply (\ref{PP}) for every $ x$ from the interval
$q_1^3\le q_2 < x\le q_3$:
$$
{\rm mes}(B_{q_2}^{[j]}\cap A^{[j]} (x)) \le {\rm
mes}(B_{q_1}^{[j]}\cap A^{[j]} (x)) \le
2^4\sigma^{[j]}_\varepsilon (x) \times {\rm mes} B_{q_1}^{[j]} \le
2^5\sigma^{[j]}_\varepsilon (x) \times {\rm mes} B_{q_2}^{[j]}
$$
(in the last inequality  we use the condition  (\ref{m1}) of Lemma
2). Now as $\frac{\log_2q_3}{\log_2 q_2}= A$ the conclusion
(\ref{m2}) of Lemma 5 in the case $j=1$ follows from Lemma 3:
$$
{\rm mes}  B_{q_3}^{[1]}\ge {\rm mes}   B_{q_2}^{[1]} \left( 1-
2^5 T_{A,\varepsilon}^{[1]} (q_2)\right) \ge
 {\rm mes}   B_{q_2}^{[1]}/2.
$$
In the case $j=2$ Lemma 5 follows from Lemma 4 by a similar
argument.

{\bf 9. Acknowledgement.} It is my great pleasure to thank Yann Bugeaud for many suggestions and comments.

\end{document}